\documentclass[12pt]{article}
\usepackage{mathrsfs}
\usepackage{amsfonts}
\usepackage{amsthm,amsmath,amssymb,anysize,longtable}

\newtheorem{theorem}{Theorem}
\newtheorem{lemma}[theorem]{Lemma}

\newtheorem{definition}[theorem]{Definition}

\usepackage{xcolor}

\usepackage[numbers,sort&compress]{natbib}
\usepackage{amsmath}
\allowdisplaybreaks[3]
\marginsize{4.2cm}{4.2cm}{3.6cm}{3cm}
\numberwithin{equation}{section}

\author{\textsc{Qingyan Ren and Liming Tang\footnote{corresponding author}}\\
\small{School of Mathematical Sciences}\\
\small{Harbin Normal University}\\
\small{150025 Harbin, China}\\
\small{E-mail: limingtang@hrbnu.edu.cn}}

\date{ }
\date{ }
\usepackage[symbol]{footmisc}

\begin{document}

\thispagestyle{empty}

\noindent{\Large
Second cohomology groups and left-symmetric algebraic structures of the generalized loop Heisenberg-Virasoro algebra}
\footnote{
The work is supported by the NSF of Hei Longjiang Province (No. LH2024A014).
}

	\bigskip
	
	 \bigskip

\begin{center}	
	{\bf
		
    Qingyan Ren\footnote{School of Mathematical Sciences, Harbin Normal University, 150025 Harbin, China;\ renqingyan@stu.hrbnu.edu.cn},
Liming Tang\footnote{School of Mathematical Sciences, Harbin Normal University, 150025 Harbin, China; \ limingtang@hrbnu.edu.cn}\footnote{corresponding author}}
\end{center}

 \begin{quotation}
{\small\noindent \textbf{Abstract}:
This is the second paper in our series of papers dedicated to the study of the generalized loop Heisenberg-Virasoro algebra. The first paper is dedicated to the study of
   maps on the generalized loop Heisenberg-Virasoro algebra, including derivations, $2$-local derivations, biderivations  the automorphism groups. The present paper is dedicated to the study of the second cohomology groups and left-symmetric algebraic structures on the generalized loop Heisenberg-Virasoro algebra. 
  We describe  the second cohomology group and left-symmetric algebra structures on the generalized loop Heisenberg-Virasoro algebra by use of the Witt algebras being Lie subalgebra of the generalized loop Heisenberg-Virasoro algebra up to isomorphism. 
  
\medskip
 \vspace{0.05cm} \noindent{\textbf{Keywords}}:
 generalized loop Heisenberg-Virasoro algebra; second cocycle; second cohomology groups; left-symmetric algebra

\medskip

\vspace{0.05cm} \noindent \textbf{Mathematics Subject Classification
2020}: 17B40, 17B65, 17B68}
\end{quotation}
 \medskip

\section*{Introduction}

The loop Heisenberg-Virasoro algebra is an important infinite dimensional Lie algebra and it plays a vital role in numerous fields of mathematics and physics. In this paper,  we mainly describe the second cohomology group of the generalized loop Heisenberg-Virasoro algebra and its left-symmetric algebraic structure.

\smallskip

As is well known, homology theory of Lie algebras plays an important role in the structure and representation theory of Lie algebras. Furthermore, the second cohomology groups can be used to construct many infinite-dimensional Lie algebras and can also compute the universal central extensions of Lie algebras. In 1998, Okovi$\acute{c}$ and Zhao K studied the second cohomology groups $H^2(W,\mathcal{F})$ of the generalized Witt algebra \cite{ref6}. In 2007, Su Y and Zhao K studied the second cohomology groups of the generalized Witt algebra (more general than that defined by Okovi$\acute{c}$ and Zhao K), then the Virasoro-type algebras are defined \cite{ref7}. In 2008, Li J, Su Y and Zhu L  determined the second cohomology groups of original deformative Schr$\ddot{o}$dinger-Virasoro algebras \cite{ref8}. In 2010, Xie W obtained the second cohomology groups of the finite-dimensional generalized Witt Lie superalgebra is trivial by computing the derivations of Witt Lie superalgebra \cite{ref9}. In 2012, Li J studied the second cohomology groups of twisted deformative Schr$\ddot{o}$dinger–Virasoro Algebras and deduced the structure and generators of the second cohomology group \cite{ref10}. In 2014,  Wu H, Wang S and Yue X described the second cohomology groups of Generalized loop Schr$\ddot{o}$dinger–Virasoro algebras \cite{ref3}. In 2018, Chen H, Fan G, Han J and Su Y studied the second cohomology groups of generalized Schr$\ddot{o}$dinger–Virasoro algebras \cite{ref12}. In 2023, Yuan L and Li J determined the second cohomology groups of $q$-deformative $W(2,2)$ algebra and obtained  $H^2(W^q,W^q)$ is two-dimensional \cite{ref13}. In 2025, Guo X, Kaygorodov I and Tang L determined the biderivations and left-symmetric algebraic structure on the Mirror Heisenberg Virasoro algebras \cite{ref20}.

\smallskip

Left-symmetric algebra also is called pre-Lie algebra. It is closely related to Lie groups and Lie algebras. In 2002, Bai C and Meng D discussed the compatible left-symmetric algebra structures on some complete Lie algebras \cite{ref14}. In 2011, Kong X, Chen H and Bai C studied the compatible left-symmetric algebra structures on $W(2,2)$ \cite{ref5}. In 2012, Chen H and Li J studied the compatible left-symmetric algebra structures on the W-algebra $W(2,2)$ with some natural grading conditions and deduced any such left-symmetric algebra contains an infinite-dimensional trivial subalgebra \cite{ref16}. In 2014, Chen H and Li J also determined the left-symmetric algebra structures on the twisted Heisenberg-Virasoro algebra \cite{ref4}. In 2022, Xu C classified all graded compatible left-symmetric algebraic structures on high rank Witt algebras and all non-graded ones satisfying a minor condition \cite{ref18}. In 2024, Chen H and Wang Q studied the compatible left-symmetric algebra structures on the finite-dimensional Witt algebra \cite{ref19}.

\smallskip

In 2024, we introduced the generalized loop Heisenberg-Virasoro algebra denoted by $\mathcal{L}(\Gamma)$ and discussed the maps on the generalized loop Heisenberg-Virasoro algebra, including derivations, $2$-local derivations, biderivations  the automorphism groups in \cite{ref2}. In this paper, inspired by \cite{ref20, ref5,ref12}, we determine the second cohomology groups of $\mathcal{L}(\Gamma)$ and its compatible left-symmetric algebra structures. The structure of this paper is as follows: Firstly, second cohomology groups on $\mathcal{L}(\Gamma)$ are determined. Then, due to the non-associativity of left-symmetric algebras, direct research on them is rather challenging, a compatible left-symmetric algebraic structure on $\mathcal{L}(\Gamma)$ with similar natural graded conditions can be determined by using the natural graded conditions of the compatible left-symmetric algebraic structures on the Witt algebra.

\section{Preliminaries}

Throughout this paper, the ground is the complex field $\mathbb{C}$, all vector spaces and algebras are over $\mathbb{C}$. $\Gamma$ denotes an abelian additive group, and $\mathbb{Z}$ denotes the set of integers. Let $\mathrm{span}\{X\}$ denote the vector space generated by the vectors in $\{X\}$. We assume that $\epsilon\in\mathbb{C}$ possesses the following properties:
\begin{eqnarray}\nonumber
\mathbf{Re}\epsilon>0,\,\epsilon^{-1}\notin\mathbb{Z}\quad or\quad\mathbf{Re}\epsilon=0,\,\mathbf{Im}\epsilon>0,
\end{eqnarray}
where $\mathbf{Re}\epsilon$ is the real part of $\epsilon$ and $\mathbf{Im}\epsilon$ is the imaginary part of $\epsilon$.

Recall that a Lie algebra is called the generalized loop Heisenberg-Virasoro algebra with basis $\left\{L_{\alpha,i},H_{\beta,j}\,|\,\alpha,\beta\in\Gamma,\,i,j\in\mathbb{Z}\right\}$, subject to the following  Lie brackets (see \cite{ref2}):
\begin{longtable}{lcl}
$[L_{\alpha,i},L_{\beta,j}]$&$=$&$(\alpha-\beta)L_{\alpha+\beta,i+j},$ \\
$[L_{\alpha,i},H_{\beta,j}]$&$=$&$-\beta H_{\alpha+\beta,i+j},$ \\
$[H_{\alpha,i},H_{\beta,j}]$&$=$&$0.$
\end{longtable}

 Let $\mathcal{LV}$ be the subalgebra of $\mathcal{L}(\Gamma)$ spanned by $\{L_{\alpha,i}\,|\,\alpha\in\Gamma,\,i\in\mathbb{Z}\}$. Then $\mathcal{LV}$ is  isomorphic to the centerless generalized loop Virasoro algebra (see \cite{ref3}). In addition, the subalgebra of $\mathcal{L}(\Gamma)$ spanned by $\{L_{\alpha,0}\,|\,\alpha\in\Gamma\}$ is isomorphic to the Witt algebra. Since $L_{0,0}$ is semisimple, then the $\Gamma$-grading on $\mathcal{L}(\Gamma)$ is given by 
$$\mathcal{L}(\Gamma)=\bigoplus_{\mu \in \Gamma}\mathcal{L}(\Gamma)_\mu,$$
where
$$\mathcal{L}(\Gamma)_\mu=\left\{x\in \mathcal{L}(\Gamma)\,|\,[L_{0,0},x]=-\mu x\right\}=\mathrm{span}_\mathbb{C}\left\{L_{\mu,j},H_{\mu,j}\,|\,j\in\mathbb{Z}\right\}.$$
Recall that a derivation $\mathcal{D}\in\mathrm{Der}\mathcal{L}(\Gamma)$ is of degree $\gamma\in\Gamma$ if $D(\mathcal{L}(\Gamma)_\alpha)\subset\mathcal{L}_{\alpha+\gamma}$ for all $\alpha\in\Gamma$. Let $(\mathrm{Der}\mathcal{L}(\Gamma))_\gamma$ be the space of all derivations of degree $\gamma$ (see \cite{ref2}).

\section{Second cohomology groups on $\mathcal{L}(\Gamma)$}
A bilinear form $\psi:\mathcal{L}\times\mathcal{L}\longrightarrow\mathbb{C}$ is called a \emph{second cocycle} on Lie algebra $\mathcal{L}$ if the following conditions are satisfied:
\begin{eqnarray}\nonumber
&&\psi(x,y)=-\psi(y,x),\\ \nonumber
&&\psi(x,[y,z])+\psi(y,[z,x])+\psi(z,[x,y])=0,\,\mbox{for all}\, x,y,z\in\mathcal{L}.
\end{eqnarray}
Denote by $C^2(\mathcal{L},\mathbb{C})$ the vector space of second cocycles on $\mathcal{L}$ (see \cite{ref7}).

For any linear function $f:\mathcal{L}\longrightarrow\mathbb{C}$, one can define a second cocycle $\psi_f$ as follows: 
\begin{eqnarray}
\psi_f(x,y)=f([x,y]),\,\mbox{for all}\,\,x,y\in\mathcal{L}.
\end{eqnarray}
Such a second cocycle is called a \emph{second coboundary} on $\mathcal{L}$. Denote by $B^2(\mathcal{L},\mathbb{C})$ the vector space of second coboundaries on $\mathcal{L}$ (see \cite{ref7}). 

The quotient space
\begin{eqnarray}\nonumber
H^2(\mathcal{L},\mathbb{C})=C^2(\mathcal{L},\mathbb{C})/B^2(\mathcal{L},\mathbb{C})
\end{eqnarray}
is called the \emph{second cohomology group} of $\mathcal{L}$ (see \cite{ref7}). 

Let $\psi$ be a second cocycle of $\mathcal{L}(\Gamma)$. A linear map $f:\mathcal{L}(\Gamma)\longrightarrow\mathbb{C}$ is defined by:
\begin{eqnarray}\nonumber
&&f(L_{\alpha,i})=
\begin{cases}
-\frac{1}{\alpha}\psi(L_{0,0},L_{\alpha,i}),\,\alpha\neq0\\ \frac{1}{2}\psi(L_{1,0},L_{-1,i}),\,\alpha=0
\end{cases}\\ \nonumber
&&f(H_{\alpha,i})=
\begin{cases}
-\frac{1}{\alpha}\psi(L_{0,0},H_{\alpha,i}),\,\alpha\neq0\\ -\psi(L_{1,0},H_{-1,i}),\,\alpha=0.
\end{cases}
\end{eqnarray}
Set $\phi=\psi-\psi_f$, where $\psi_f$ is defined in (2.1). For all $x,y\in\mathcal{L}(\Gamma)$, we have $\phi(x,y)=-\phi(y,x),\,x,y\in\mathcal{L}(\Gamma)$.

Furthermore, we have
\begin{eqnarray}
&&\phi(L_{0,0},L_{\alpha,i})=0\,\,\mbox{for all}\,\,0\neq\alpha\in\Gamma,\\ 
&&\phi(L_{1,0},L_{-1,i})=0\,\,\mbox{for all}\,\, i\in\mathbb{Z},\\ 
&&\phi(L_{0,0},H_{\alpha,i})=0\,\,\mbox{for all}\,\,0\neq\alpha\in\Gamma,\\ 
&&\phi(L_{1,0},H_{-1,i})=0\,\,\mbox{for all}\,\,i\in\mathbb{Z}.
\end{eqnarray}

The second cohomology group on $\mathcal{L}(\Gamma)$ is determined by the following theorem. 

\begin{theorem}
The \emph{second cohomology group} on $\mathcal{L}(\Gamma)$
$$H^2(\mathcal{L}(\Gamma),\mathbb{C})=\prod_{\substack{k\in\mathbb{Z}\\x\in X}}\mathbb{C}\overline{\phi}_{k,x}$$
is the direct product of all $\mathbb{C}\overline{\phi}_{k,x}$, where $\overline{\phi}_{k,x}$ is the cohomology class of $\phi_{k,x}$ defined by
\begin{eqnarray}\nonumber
\phi_{k,1}(L_{\alpha,i},L_{\beta,j})&=&\delta_{\alpha+\beta,0}\delta_{i+j,k}\frac{\alpha^3-\alpha}{12},\\ \nonumber
\phi_{k,2}(L_{\alpha,i},H_{\beta,j})&=&\delta_{\alpha+\beta,0}\delta_{i+j,k}(\alpha^2-\alpha),\\ \nonumber
\phi_{k,3}(H_{\alpha,i},H_{\beta,j})&=&\delta_{\alpha+\beta,0}\delta_{i+j,k}\alpha,
\end{eqnarray}
where $k\in\mathbb{Z}$, $x\in X=\{1,2,3\}$ and $\delta_{i,j}$ denotes $1$ when $i=j$ and $0$ when $i\neq j$. 
\end{theorem}

\subsection{Proof of Theorem 1}
For simplicity, denote
\begin{eqnarray}\nonumber
&&A_{\alpha,i,j}=\phi(L_{\alpha,i},L_{-\alpha,j}),\\ \nonumber
&&B_{\alpha,i,j}=\phi(L_{\alpha,i},H_{-\alpha,j}),\\ \nonumber
&&C_{\alpha,i,j}=\phi(H_{\alpha,i},H_{-\alpha,j}),\\ \nonumber
&&A_{\alpha,\beta,i,j}=\phi(L_{\alpha,i},L_{\beta,j}),\\ \nonumber
&&B_{\alpha,\beta,i,j}=\phi(L_{\alpha,i},H_{\beta,j}),\\ \nonumber
&&C_{\alpha,\beta,i,j}=\phi(H_{\alpha,i},H_{\beta,j}),\,\mbox{for all}\,\alpha,\beta\in\Gamma,\,i,j\in\mathbb{Z}.
\end{eqnarray}
\begin{lemma}
$A_{\alpha,i,j+k}=A_{\alpha,j,i+k}, \,\mbox{for all}\,\,i,j,k\in\mathbb{Z},\,0\neq \alpha\in\Gamma$.

\begin{proof}
Since \begin{eqnarray}\nonumber
(\alpha-\beta)A_{\alpha+\beta,i+j,k}&=&(\alpha-\beta)\phi(L_{\alpha+\beta,i+j},L_{-(\alpha+\beta),k})\\ \nonumber
&=&\phi([L_{\alpha,i},L_{\beta,j}],L_{-(\alpha+\beta),k})\\ \nonumber
&=&\phi(L_{\alpha,i},[L_{\beta,j},L_{-(\alpha+\beta),k}])+\phi(L_{\beta,j},[L_{-(\alpha+\beta),k},L_{\alpha,i}])\\ \nonumber
&=&(2\beta+\alpha)\phi(L_{\alpha,i},L_{-\alpha,j+k})-(2\alpha+\beta)\phi(L_{\beta,j},L_{-\beta,i+k})\\ \nonumber
&=&(2\beta+\alpha)A_{\alpha,i,j+k}-(2\alpha+\beta)A_{\beta,j,i+k}
\end{eqnarray}
Let $\alpha=\beta\neq0$. Then we have
$$A_{\alpha,i,j+k}=A_{\alpha,j,i+k}.$$
\end{proof}
\end{lemma}

\begin{lemma}
$A_{0,i,j}=0, \,\mbox{for all}\,i,j\in\mathbb{Z}$.\

\begin{proof}
Since\begin{eqnarray}\nonumber
2A_{0,i,j}&=&2\phi(L_{0,i},L_{0,j})\\ \nonumber
&=&\phi(L_{0,i},[L_{1,0},L_{-1,j}])\\ \nonumber
&=&\phi([L_{-1,j},L_{0,i}],L_{1,0})+\phi([L_{0,i},L_{1,0}],L_{-1,j})\\ \nonumber
&=&A_{1,0,i+j}-A_{1,i,j}
\end{eqnarray}
By Lemma 2, we have $A_{1,0,i+j}=A_{1,i,j}.$
Then
$$A_{0,i,j}=0, \,\mbox{for all}\,\,i,j\in\mathbb{Z}.$$
\end{proof}
\end{lemma}

Since  $A_{\alpha,i,j}$ only depends on the sum $i+j$. Then in the sequel we shortly write $A_{\alpha,i,j}$ as $A_{\alpha,i+j}$.

\begin{lemma}
$A_{\alpha,\beta,i,j}=\delta_{\alpha+\beta,0}\frac{\alpha^3-\alpha}{6}A_{2,i+j},\,\mbox{for all}\,\,\alpha,\beta\in\Gamma,\,i,j\in\mathbb{Z}$.

\begin{proof}Since
\begin{eqnarray}\nonumber
-\beta A_{\alpha,\beta,i,j}&=&-\beta\phi(L_{\alpha,i},L_{\beta,j})\\ \nonumber
&=&\phi(L_{\alpha,i},[L_{0,0},L_{\beta,j}])\\ \nonumber
&=&\phi([L_{\alpha,i},L_{0,0}],L_{\beta,j})+\phi(L_{0,0},[L_{\alpha,i},L_{\beta,j}])\\ \nonumber
&=&\alpha\phi(L_{\alpha,i},L_{\beta,j})+(\alpha-\beta)\phi(L_{0,0},L_{\alpha+\beta,i+j})\\ \nonumber
&=&\alpha A_{\alpha,\beta,i,j}
\end{eqnarray}
By (2.2), then we have $\phi(L_{0,0},L_{\alpha+\beta,i+j})=0,\,\alpha+\beta\neq 0$. Thus $A_{\alpha,\beta,i,j}=0,\,\alpha+\beta\neq0$.
We obtain
$$\phi([L_{\alpha,i},L_{\beta,j}],L_{-(\alpha+\beta),k})+\phi([L_{\beta,j},L_{-(\alpha+\beta),k}],L_{\alpha,i})+\phi([L_{-(\alpha+\beta),k},L_{\alpha,i}],L_{\beta,j})=0.$$
Simplifying the above formula gives that
\begin{eqnarray}
(\alpha-\beta)A_{\alpha+\beta,i+j,k}-(\alpha+2\beta)A_{\alpha,i,j+k}+(2\alpha+\beta)A_{\beta,j,i+k}=0.
\end{eqnarray}
Taking $\beta=1$ in (2.6). By (2.3) and Lemma 3, $A_{1,j,i+k}=0$. We have
$$(\alpha-1)A_{\alpha+1,i+j+k}-(\alpha+2)A_{\alpha,i+j+k}=0.$$
Replace $\alpha$ with $\alpha-1$ in the above formula. We have
$$(\alpha-2)A_{\alpha,i+j+k}-(\alpha+1)A_{\alpha-1,i+j+k}=0.$$
Taking $\alpha=\alpha-1,\,\beta=2$ in (2.6). We have
$$(\alpha-3)A_{\alpha+1,i+j+k}-(\alpha+3)A_{\alpha-1,i+j+k}=-2\alpha A_{2,i+j+k}.$$
We get
$$A_{\alpha,i+j+k}=\frac{\alpha^3-\alpha}{6}A_{2,i+j+k}.$$
Then
$$A_{\alpha,\beta,i,j}=\delta_{\alpha+\beta,0}\frac{\alpha^3-\alpha}{6}A_{2,i+j}.$$
\end{proof}
\end{lemma}

In the following, we calculate $\phi_{k,2}(L_{\alpha,i},H_{\beta,j})$.

\begin{lemma}
$B_{\alpha,i,j+k}=B_{\alpha,j,i+k},\,\mbox{for all}\,\,i,j,k\in\mathbb{Z},\,\alpha\in\Gamma,\,\alpha\neq0$.

\begin{proof}
Since \begin{eqnarray}\nonumber
(\alpha-\beta)B_{\alpha+\beta,i+j,k}&=&(\alpha-\beta)\phi(L_{\alpha+\beta,i+j},H_{-(\alpha+\beta),k})\\ \nonumber
&=&\phi([L_{\alpha,i},L_{\beta,j}],H_{-(\alpha+\beta),k})\\ \nonumber
&=&\phi(L_{\alpha,i},[L_{\beta,j}, H_{-(\alpha+\beta),k}])+\phi(L_{\beta,j},[H_{-(\alpha+\beta),k},L_{\alpha,i}])\\ \nonumber
&=&(\alpha+\beta)B_{\alpha,i,j+k}-(\alpha+\beta)B_{\beta,j,i+k}
\end{eqnarray}
Let $\alpha=\beta$. Then we have
$$0=2\alpha B_{\alpha,i,j+k}-2\alpha B_{\alpha,j,i+k}.$$
Since $\alpha\neq 0$, then
$$B_{\alpha,i,j+k}=B_{\alpha,j,i+k}.$$
\end{proof}
\end{lemma}

\begin{lemma}
$B_{0,i,j}=0,\,\mbox{for all}\,\,i,j\in\mathbb{Z}$.

\begin{proof}
Since\begin{eqnarray}\nonumber
B_{0,i,j}&=&\phi(L_{0,i},H_{0,j})\\ \nonumber
&=&\phi(L_{0,i},[L_{1,0},H_{-1,j}])\\ \nonumber
&=&\phi([H_{-1,j},L_{0,i}],L_{1,0})+\phi([L_{0,i},L_{1,0}],H_{-1,j})\\ \nonumber
&=&B_{1,0,i+j}-B_{1,i,j}
\end{eqnarray}
By (2.5), then we have $B_{1,0,i+j=0}=0$. By Lemma 5, we obtain $B_{1,i,j}=B_{1,0,i+j}=0$, then we conclude that $B_{0,i,j}=0$.
\end{proof}
\end{lemma}

It follows from the above two lemmas that $B_{\alpha,i,j}$ only depends on the sum $i+j$. Thus in the sequel we shortly write $B_{\alpha,i,j}$ as $B_{\alpha,i+j}$.

\begin{lemma}
$B_{\alpha,\beta,i,j}=\delta_{\alpha+\beta,0}\frac{\alpha^2-\alpha}{2}B_{2,i+j},\,\mbox{for all}\,\,\alpha,\beta\in\Gamma,\,i,j\in\mathbb{Z}$.

\begin{proof}
Since \begin{eqnarray}\nonumber
-\beta B_{\alpha,\beta,i,j}&=&-\beta\phi(L_{\alpha,i},H_{\beta,j})\\ \nonumber
&=&\phi(L_{\alpha,i},[L_{0,0},H_{\beta,j}])\\ \nonumber
&=&\phi([L_{\alpha,i},L_{0,0}],H_{\beta,j})+\phi(L_{0,0},[L_{\alpha,i},H_{\beta,j}])\\ \nonumber
&=&\alpha B_{\alpha,\beta,i,j}-\beta B_{0,\alpha+\beta,0,i+j}
\end{eqnarray}
By (2.4), then we have $B_{0,\alpha+\beta,0,i+j}$. Thus
$$B_{\alpha,\beta,i,j}=0,\,\alpha+\beta\neq0.$$
We obtain
$$\phi([L_{\alpha,i},L_{\beta,j}],H_{-(\alpha+\beta),k})+\phi([L_{\beta,j},H_{-(\alpha+\beta),k}],L_{\alpha,i})+\phi([H_{-(\alpha+\beta),k},L_{\alpha,i}],L_{\beta,j})=0.$$
Simplifying the above formula gives that
\begin{eqnarray}
(\alpha-\beta)B_{\alpha+\beta,i+j,k}-(\alpha+\beta)B_{\alpha,i,j+k}+(\alpha+\beta)B_{\beta,j,i+k}=0.
\end{eqnarray}
Taking $\beta=1$ in (2.7). By (2.5) and Lemma 6, $B_{1,j,i+k}=0$. We have
$$(\alpha-1)B_{\alpha+1,i+j+k}-(\alpha+1)B_{\alpha,i+j+k}=0.$$
Replace $\alpha$ with $\alpha-1$ in the above formula,
$$(\alpha-2)B_{\alpha,i+j+k}-\alpha B_{\alpha-1,i+j+k}=0.$$
Taking $\alpha=\alpha-1,\,\beta=2$ in (2.7). We have
$$(\alpha-3)B_{\alpha+1,i+j+k}-(\alpha+1)B_{\alpha-1,i+j+k}+(\alpha+1)B_{2,i+j+k}=0.$$
We get
$$B_{\alpha,i+j+k}=\frac{\alpha^2-\alpha}{2}B_{2,i+j+k},\,\alpha\neq1.$$
Taking $\alpha=2,\,\beta=-1$ in (2.7), we have
$$3B_{1,i+j,k}-B_{2,i,j+k}+B_{-1,j,i+k}=0.$$
By (2.5), we get $B_{-1,i+j+k}=B_{-1,j,i+k}=B_{2,i,j+k}=B_{2,i+j+k}.$
We have 
$$B_{\alpha,\beta,i,j}=B_{2,i+j},\,\alpha=-1.$$
Thus
$$B_{\alpha,\beta,i,j}=\delta_{\alpha+\beta,0}\frac{\alpha^2-\alpha}{2}B_{2,i+j}.$$
\end{proof}
\end{lemma}

Finally, we calculate $\phi_{k,3}(H_{\alpha,i},H_{\beta,j})$.

\begin{lemma}
$C_{\alpha,i,j+k}=C_{\alpha,j,i+k},\,\mbox{for all}\,\,i,j,k\in\mathbb{Z},\,\alpha\in\Gamma,\,\alpha\neq0$.

\begin{proof}
Since \begin{eqnarray}\nonumber
-\alpha C_{\alpha+\beta,i+j,k}&=&-\alpha\phi(H_{\alpha+\beta,i+j},H_{-(\alpha+\beta),k})\\ \nonumber
&=&\phi([L_{\beta,i},H_{\alpha,j}],H_{-(\alpha+\beta),k})\\ \nonumber
&=&\phi(L_{\beta,i},[H_{\alpha,i},H_{-(\alpha+\beta),k}])+\phi(H_{\alpha,j},[H_{-(\alpha+\beta),k},L_{\beta,i}])\\ \nonumber
&=&-(\alpha+\beta)C_{\alpha,j,i+k}
\end{eqnarray}
Let $\beta=0$. Then we have
$$C_{\alpha,i+j,k}=C_{\alpha,j,i+k}.$$
Similarly, we have
$$C_{\alpha,i+j,k}=C_{\alpha,i,j+k}.$$
Then
$$C_{\alpha,i,j+k}=C_{\alpha,j,i+k}.$$
\end{proof}
\end{lemma}

\begin{lemma}
$C_{0,i,j}=0,\,\mbox{for all}\,\,i,j\in\mathbb{Z}$

\begin{proof}
A direct computation shows \begin{eqnarray}\nonumber
C_{0,i,j}&=&\phi(H_{0,i},H_{0,j})\\ \nonumber
&=&\phi(H_{0,i},[L_{1,0},H_{-1,j}])\\ \nonumber
&=&\phi([H_{-1,j},H_{0,j}],L_{1,0})+\phi([H_{0,i},L_{1,0}],H_{-1,j})\\ \nonumber
&=&0
\end{eqnarray}
\end{proof}
\end{lemma}

It follows from the above two lemmas that $C_{\alpha,i,j}$ only depends on the sum $i+j$. Then in the sequel we shortly write $C_{\alpha,i,j}$ as $C_{\alpha,i+j}$.

\begin{lemma}
$C_{\alpha,\beta,i,j}=\delta_{\alpha+\beta,0}\alpha C_{1,i+j+k},\,\mbox{for all}\,\,\alpha,\beta\in\Gamma,\,i,j\in\mathbb{Z}$

\begin{proof}
Since\begin{eqnarray}\nonumber
-\alpha C_{\alpha,\beta,i,j}&=&-\alpha\phi(H_{\alpha,i},H_{\beta,j})\\ \nonumber
&=&\phi([L_{0,0},H_{\alpha,i}],H_{\beta,j})\\ \nonumber
&=&\phi(L_{0,0},[H_{\alpha,i},H_{\beta,j}])+\phi(H_{\alpha,i},[H_{\beta,j},L_{0,0}])\\ \nonumber
&=&\beta C_{\alpha,\beta,i,j}
\end{eqnarray}
Then
$$C_{\alpha,\beta,i,j}=0,\,\alpha+\beta\neq0.$$
We obtain
$$\phi([L_{\alpha,i},H_{\beta,j}],H_{-(\alpha+\beta),k})+\phi([H_{\beta,j},H_{-(\alpha+\beta),k}],L_{\alpha,i})+\phi([H_{-(\alpha+\beta),k},L_{\alpha,i}],H_{\beta,j})=0.$$
Simplifying the above formula gives that
$$-\beta C_{\alpha+\beta,i+j+k}+(\alpha+\beta)C_{\beta,i+j+k}=0.$$
Let $\alpha=\alpha-1,\,\beta=1$. Then we have
$$C_{\alpha,i+j+k}=\alpha C_{1,i+j+k}.$$
Thus
$$C_{\alpha,\beta,i,j}=\delta_{\alpha+\beta,0}\alpha C_{1,i+j+k}.$$
\end{proof}
\end{lemma}

Now, it follows from \cite{ref3} and Lemmas 2-10.
For any $k\in\mathbb{C}$, define
\begin{eqnarray}\nonumber
\phi_{k,1}(L_{\alpha,i},L_{\beta,j})&=&\delta_{\alpha+\beta,0}\delta_{i+j,k}\frac{\alpha^3-\alpha}{12},\\ \nonumber
\phi_{k,2}(L_{\alpha,i},H_{\beta,j})&=&\delta_{\alpha+\beta,0}\delta_{i+j,k}(\alpha^2-\alpha),\\ \nonumber
\phi_{k,3}(H_{\alpha,i},H_{\beta,j})&=&\delta_{\alpha+\beta,0}\delta_{i+j,k}\alpha.
\end{eqnarray}
Then we have
$$H^2(\mathcal{L}(\Gamma),\mathbb{C})=\prod_{\substack{k\in\mathbb{Z}\\x\in X}}\mathbb{C}\overline{\phi}_{k,x}.$$

\section{Left-symmetric algebraic structures on $\mathcal{L}(\Gamma)$}
\begin{definition}(see \cite{ref4})
Let $\mathcal{A}$ be a vector space over a field $\mathbb{F}$ equipped with a bilinear product $(x,y)\mapsto xy$, $\mathcal{A}$ is called a \emph{left-symmetric algebra} if any $x,y,z\in\mathcal{A}$, the associator
$$(x,y,z)=(xy)z-x(yz)$$
is symmetric in $x,y$, that is,
$$(x,y,z)=(y,x,z)\,\mbox{or equivalently}\,\,(xy)z-x(yz)=(yx)z-y(xz).$$
\end{definition}

Recall that $\mathcal{A}$ is a left-symmetric algebra. For any $x\in\mathcal{A}$, denote by $L_x$ the left multiplication operator (i.e. $L_x(y)=xy,\,y\in\mathcal{A}$). The commutator
\begin{eqnarray}\nonumber
[x,y]=xy-yx,\,\mbox{for all}\, x,y\in\mathcal{A},
\end{eqnarray}
defines a Lie algebra $\mathcal{G}(\mathcal{A})$, which is called the sub-adjacent Lie algebra of $\mathcal{A}$ and $\mathcal{A}$ is called a compatible left-symmetric algebra structure on $\mathcal{G}(\mathcal{A})$ (see \cite{ref4}).

Let $\mathcal{G}$ be a Lie algebra. Denote $\mathcal{A}(\mathcal{G})$ by the compatible left-symmetric algebra of $\mathcal{G}$. The Witt algebra $\mathcal{W}=\mathrm{span}\{L_{n,0}\,|\,n\in\mathbb{Z}\}$ is said to have the natural grading condition if nonzero multiplications of basis element of $\mathcal{A}(\mathcal{W})$ satisfies (see \cite{ref5})
\begin{eqnarray}
L_{m,0}L_{n,0}=f(m,n)L_{m+n,0},\,m,n\in\mathbb{Z}.
\end{eqnarray}

\begin{lemma}(see \cite{ref5})
Any compatible left-symmetric algebra structure on the Witt algebra satisfying (3.1) is given by the multiplication table:
\begin{eqnarray}\nonumber
L_{m,0}L_{n,0}=\frac{\alpha+n+\alpha\epsilon m}{1+\epsilon(m+n)}L_{m+n,0},
\end{eqnarray}
where $\alpha,\epsilon\in\mathbb{C}$ and $\mathbf{Re}\epsilon>0,\,\epsilon^{-1}\notin\mathbb{Z}$ or $\mathbf{Re}\epsilon=0,\,\mathbf{Im}\epsilon>0$.
\end{lemma}

The Witt algera $\mathcal{W}$ is a subalgebra of the generalized Loop Heisenberg-Virasoro algebra $\mathcal{L}(\Gamma)$. Naturally, it can be assumed that the left-symmetric algebra structure of $\mathcal{L}(\Gamma)$ satisfies the similar relation in (3.1). Thus we can denote as following: 
\begin{eqnarray}\nonumber
L_{\alpha,i}L_{\beta,j}&=&a(\alpha,i,\beta,j)L_{\alpha+\beta,i+j},\\ \nonumber
L_{\alpha,i}H_{\beta,j}&=&b(\alpha,i,\beta,j)H_{\alpha+\beta,i+j},\\
H_{\alpha,i}L_{\beta,j}&=&c(\alpha,i,\beta,j)H_{\alpha+\beta,i+j},\\ \nonumber
H_{\alpha,i}H_{\beta,j}&=&d(\alpha,i,\beta,j)L_{\alpha+\beta,i+j}+e(\alpha,i,\beta,j)H_{\alpha+\beta,i+j},\,\alpha,\beta\in\Gamma,\,i,j\in\mathbb{Z},
\end{eqnarray}
where $a(\alpha,i,\beta,j),b(\alpha,i,\beta,j),c(\alpha,i,\beta,j),d(\alpha,i,\beta,j),e(\alpha,i,\beta,j)$ are the complex-valued functions.

\begin{theorem}
Any compatible left-symmetric algebraic structure on the generalized loop Heisenberg-Virasoro algebra satisfying the relation (3.3) is given by the following functions:
\begin{eqnarray}\nonumber
&&a(\alpha,i,\beta,j)=\frac{-\beta(1+\epsilon\beta)}{1+\epsilon(\alpha+\beta)},\\ \nonumber
&&b(\alpha,i,\beta,j)=-\beta(1+(1-\epsilon\beta)m\delta_{\alpha+\beta,0}),\\ \nonumber
&&c(\alpha,i,\beta,j)=\beta(1+\epsilon\beta)m\delta_{\alpha+\beta,0},\\ \nonumber
&&d(\alpha,i,\beta,j)=0,\\ \nonumber
&&e(\alpha,i,\beta,j)=0,
\end{eqnarray}
where $\alpha,\beta\in\Gamma,\,i,j\in\mathbb{Z}$ and $\mathbf{Re}\epsilon>0,\,\epsilon^{-1}\notin\mathbb{Z}$ or $\mathbf{Re}\epsilon=0,\,\mathbf{Im}\epsilon>0$.
\end{theorem}
The algebra obtained in the previous theorem is called the graded generalized loop Heisenberg-Virasoro left-symmetric algebra $\mathcal{A}$.

\subsection{Proof of Theorem 2}
\begin{lemma}
Let $\mathcal{L}(\Gamma)$ be the generalized loop Heisenberg-Virasoro algebra. A bilinear product induced from by (3.2) gives a compatible left-symmetric algebra structure on $\mathcal{L}(\Gamma)$ if and only if 
\begin{eqnarray}\nonumber
&&a(\alpha,i,\beta,j)-a(\beta,j,\alpha,i)=\alpha-\beta,\\ \nonumber
&&b(\alpha,i,\beta,j)-c(\beta,j,\alpha,i)=-\beta,\\ \nonumber
&&d(\alpha,i,\beta,j)=d(\beta,j,\alpha,i),\,e(\alpha,i,\beta,j)=e(\beta,j,\alpha,i),\\ \nonumber
&&a(\beta,j,\gamma,k)a(\alpha,i,\beta+\gamma,j+k)-a(\alpha,i,\gamma,k)a(\beta,j,\alpha+\gamma,i+k)=(\alpha-\beta)a(\alpha+\beta,i+j,\gamma,k),\\ \nonumber
&&b(\beta,j,\gamma,k)b(\alpha,i,\beta+\gamma,j+k)-a(\alpha,i,\gamma,k)a(\beta,j,\alpha+\gamma,i+k)=(\alpha-\beta)b(\alpha+\beta,i+j,\gamma,k),\\ \nonumber
&&c(\beta,j,\gamma,k)b(\alpha,i,\beta+\gamma,j+k)-a(\alpha,i,\gamma,k)c(\beta,j,\alpha+\gamma,i+k)=-\beta c(\alpha+\beta,i+j,\gamma,k),\\ \nonumber
&&d(\beta,j,\gamma,k)a(\alpha,i,\beta+\gamma,j+k)-b(\alpha,i,\gamma,k)d(\beta,j,\alpha+\gamma,i+k)=-\beta d(\alpha+\beta,i+j,\gamma,k),\\ \nonumber
&&e(\beta,j,\gamma,k)b(\alpha,i,\beta+\gamma,j+k)-b(\alpha,i,\gamma,k)e(\beta,j,\alpha+\gamma,i+k)=-\beta e(\alpha+\beta,i+j,\gamma,k),\\ 
&&c(\beta,j,\gamma,k)d(\alpha,i,\beta+\gamma,j+k)-c(\alpha,i,\gamma,k)d(\beta,j,\alpha+\gamma,i+k)=0,\\ \nonumber
&&c(\beta,j,\gamma,k)e(\alpha,i,\beta+\gamma,j+k)-c(\alpha,i,\gamma,k)e(\beta,j,\alpha+\gamma,i+k)=0,\\ \nonumber
&&e(\beta,j,\gamma,k)e(\alpha,i,\beta+\gamma,j+k)+d(\alpha,i,\gamma,k)c(\beta,j,\alpha+\gamma,i+k)=d(\beta,j,\gamma,k)c(\alpha,i,\beta+\gamma,j+k)\\ \nonumber
&&+e(\alpha,i,\gamma,k)e(\beta,j,\alpha+\gamma,i+k),\\ \nonumber
&&e(\beta,j,\gamma,k)d(\alpha,i,\beta+\gamma,j+k)-e(\alpha,i,\gamma,k)d(\beta,j,\alpha+\gamma,i+k)=0.
\end{eqnarray}

\begin{proof}
(\romannumeral1) By
\begin{eqnarray}\nonumber  
[L_{\alpha,i},L_{\beta,j}]&=&L_{\alpha,i}L_{\beta,j}-L_{\beta,j}L_{\alpha,i}\\ \nonumber
&=&a(\alpha,i,\beta,j)L_{\alpha+\beta,i+j}-a(\beta,j,\alpha,i)L_{\alpha+\beta,i+j},\\ \nonumber
[L_{\alpha,i},L_{\beta,j}]&=&(\alpha-\beta)L_{\alpha+\beta,i+j},
\end{eqnarray}
we can conclude that $$a(\alpha,i,\beta,j)-a(\beta,j,\alpha,i)=\alpha-\beta.$$

(\romannumeral2) By
\begin{eqnarray}\nonumber  
[L_{\alpha,i},H_{\beta,j}]&=&L_{\alpha,i}H_{\beta,j}-H_{\beta,j}L_{\alpha,i}\\ \nonumber
&=&b(\alpha,i,\beta,j)H_{\alpha+\beta,i+j}-c(\beta,j,\alpha,i)H_{\alpha+\beta,i+j},\\ \nonumber
[L_{\alpha,i},H_{\beta,j}]&=&-\beta H_{\alpha+\beta,i+j},
\end{eqnarray}
we can conclude that$$b(\alpha,i,\beta,j)-c(\beta,j,\alpha,i)=-\beta.$$

(\romannumeral3) By
\begin{eqnarray}\nonumber  
[H_{\alpha,i},H_{\beta,j}]&=&H_{\alpha,i}H_{\beta,j}-H_{\beta,j}H_{\alpha,i}\\ \nonumber
&=&d(\alpha,i,\beta,j)L_{\alpha+\beta,i+j}+e(\alpha,i,\beta,j)H_{\alpha+\beta,i+j}\\ \nonumber
&&-d(\beta,j,\alpha,i)L_{\alpha+\beta,i+j}-e(\beta,j,\alpha,i)H_{\alpha+\beta,i+j},\\ \nonumber
[H_{\alpha,i},H_{\beta,j}]&=&0,
\end{eqnarray}
we can conclude that
$$d(\alpha,i,\beta,j)=d(\beta,j,\alpha,i),\,e(\alpha,i,\beta,j)=e(\beta,j,\alpha,i).$$

(\romannumeral4) By 
$$(L_{\alpha,i},L_{\beta,j},L_{\gamma,k})=(L_{\beta,j},L_{\alpha,i},L_{\gamma,k}),$$
we have
$$(L_{\alpha,i}L_{\beta,j})L_{\gamma,k}-L_{\alpha,i}(L_{\beta,j}L_{\gamma,k})=(L_{\beta,j}L_{\alpha,i})L_{\gamma,k}-L_{\beta,j}(L_{\alpha,i}L_{\gamma,k}).$$
Hence,
\begin{eqnarray}\nonumber
&&a(\alpha,i,\beta,j)L_{\alpha+\beta,i+j}L_{\gamma,k}-a(\beta,j,\gamma,k)L_{\alpha,i}L_{\beta+\gamma,j+k}\\ \nonumber
&=&a(\beta,j,\alpha,i)L_{\alpha+\beta,i+j}L_{\gamma,k}-a(\alpha,i,\gamma,k)L_{\beta,j}L_{\alpha+\gamma,i+k}.
\end{eqnarray}
Then we have
\begin{eqnarray}\nonumber
&&a(\alpha,i,\beta,j)a(\alpha+\beta,i+j,\gamma,k)-a(\beta,j,\gamma,k)a(\alpha,i,\beta+\gamma,j+k)\\ \nonumber
&=&a(\beta,j,\alpha,i)a(\alpha+\beta,i+j,\gamma,k)-a(\alpha,i,\gamma,k)a(\beta,j,\alpha+\gamma,i+k).
\end{eqnarray}
By (\romannumeral1), we can conclude that
$$a(\beta,j,\gamma,k)a(\alpha,i,\beta+\gamma,j+k)-a(\alpha,i,\gamma,k)a(\beta,j,\alpha+\gamma,i+k)=(\alpha-\beta)a(\alpha+\beta,i+j,\gamma,k).$$

(\romannumeral5) By $$(L_{\alpha,i},L_{\beta,j},H_{\gamma,k})=(L_{\beta,j},L_{\alpha,i},H_{\gamma,k}),$$
we have
$$(L_{\alpha,i}L_{\beta,j})H_{\gamma,k}-L_{\alpha,i}(L_{\beta,j}H_{\gamma,k})=(L_{\beta,j}L_{\alpha,i})H_{\gamma,k}-L_{\beta,j}(L_{\alpha,i}H_{\gamma,k}).$$
Hence
\begin{eqnarray}\nonumber
&&a(\alpha,i,\beta,j)L_{\alpha+\beta,i+j}H_{\gamma,k}-b(\beta,j,\gamma,k)L_{\alpha,i}H_{\beta+\gamma,j+k}\\ \nonumber
&=&a(\beta,j,\alpha,i)L_{\alpha+\beta,i+j}H_{\gamma,k}-b(\alpha,i,\gamma,k)L_{\beta,j}H_{\alpha+\gamma,i+k}.
\end{eqnarray}
Then we have
\begin{eqnarray}\nonumber
&&a(\alpha,i,\beta,j)b(\alpha+\beta,i+j,\gamma,k)-b(\beta,j,\gamma,k)b(\alpha,i,\beta+\gamma,j+k)\\ \nonumber
&=&a(\beta,j,\alpha,i)b(\alpha+\beta,i+j,\gamma,k)-b(\alpha,i,\gamma,k)b(\beta,j,\alpha+\gamma,i+k).
\end{eqnarray}
By (\romannumeral2), we can conclude that
$$b(\beta,j,\gamma,k)b(\alpha,i,\beta+\gamma,j+k)-b(\alpha,i,\gamma,k)b(\beta,j,\alpha+\gamma,i+k)=(\alpha-\beta)b(\alpha+\beta,i+j,\gamma,k).$$

(\romannumeral6) By 
$$(L_{\alpha,i},H_{\beta,j},L_{\gamma,k})=(H_{\beta,j},L_{\alpha,i},L_{\gamma,k}),$$
we have
$$(L_{\alpha,i}H_{\beta,j})L_{\gamma,k}-L_{\alpha,i}(H_{\beta,j}L_{\gamma,k})=(H_{\beta,j}L_{\alpha,i})L_{\gamma,k}-H_{\beta,j}(L_{\alpha,i}L_{\gamma,k}),$$
Hence
\begin{eqnarray}\nonumber
&&b(\alpha,i,\beta,j)H_{\alpha+\beta,i+j}L_{\gamma,k}-c(\beta,j,\gamma,k)L_{\alpha,i}H_{\beta+\gamma,j+k}\\ \nonumber
&=&c(\beta,j,\alpha,i)H_{\alpha+\beta,i+j}L_{\gamma,k}-a(\alpha,i,\gamma,k)H_{\beta,j}L_{\alpha+\gamma,i+k}.
\end{eqnarray}
Then we have
\begin{eqnarray}\nonumber
&&b(\alpha,i,\beta,j)c(\alpha+\beta,i+j,\gamma,k)-c(\beta,j,\gamma,k)b(\alpha,i,\beta+\gamma,j+k)\\ \nonumber
&=&c(\beta,j,\alpha,i)c(\alpha+\beta,i+j,\gamma,k)-a(\alpha,i,\gamma,k)c(\beta,j,\alpha+\gamma,i+k).
\end{eqnarray}
By (\romannumeral3), we can conclude that
$$c(\beta,j,\gamma,k)b(\alpha,i,\beta+\gamma,j+k)-a(\alpha,i,\gamma,k)c(\beta,j,\alpha+\gamma,i+k)=-\beta c(\alpha+\beta,i+j,\gamma,k).$$

(\romannumeral7) By
$$(L_{\alpha,i},H_{\beta,j},H_{\gamma,k})=(H_{\beta,j},L_{\alpha,i},H_{\gamma,k}),$$
we have
$$(L_{\alpha,i}H_{\beta,j})H_{\gamma,k}-L_{\alpha,i}(H_{\beta,j}H_{\gamma,k})=(H_{\beta,j}L_{\alpha,i})H_{\gamma,k}-H_{\beta,j}(L_{\alpha,i}H_{\gamma,k}).$$
Hence
\begin{eqnarray}\nonumber
&&b(\alpha,i,\beta,j)H_{\alpha+\beta,i+j}H_{\gamma,k}-d(\beta,j,\gamma,k)L_{\alpha,i}L_{\beta+\gamma,j+k}-e(\beta,j,\gamma,k)L_{\alpha,i}H_{\beta+\gamma,j+k}\\ \nonumber
&=&c(\beta,j,\alpha,i)H_{\alpha+\beta,i+j}H_{\gamma,k}-b(\alpha,i,\gamma,k)H_{\beta,j}H_{\alpha+\gamma,i+k}.
\end{eqnarray}
Then we have
\begin{eqnarray}\nonumber
&&d(\alpha+\beta,i+j,\gamma,k)(b(\alpha,i,\beta,j)-c(\beta,j,\alpha,i))\\ \nonumber
&=&d(\beta,j,\gamma,k)a(\alpha,i,\beta+\gamma,j+k)-b(\alpha,i,\gamma,k)d(\beta,j,\alpha+\gamma,j+k),\\ \nonumber
&&e(\alpha+\beta,i+j,\gamma,k)(b(\alpha,i,\beta,j)-c(\beta,j,\alpha,i))\\ \nonumber
&=&e(\beta,j,\gamma,k)b(\alpha,i,\beta+\gamma,j+k)-b(\alpha,i,\gamma,k)e(\beta,j,\alpha+\gamma,j+k).
\end{eqnarray}
By (\romannumeral2), we can conclude that
\begin{eqnarray}\nonumber
&&d(\beta,j,\gamma,k)a(\alpha,i,\beta+\gamma,j+k)-b(\alpha,i,\gamma,k)d(\beta,j,\alpha+\gamma,i+k)=-\beta d(\alpha+\beta,i+j,\gamma,k),\\ \nonumber
&&e(\beta,j,\gamma,k)b(\alpha,i,\beta+\gamma,j+k)-b(\alpha,i,\gamma,k)e(\beta,j,\alpha+\gamma,i+k)=-\beta e(\alpha+\beta,i+j,\gamma,k).
\end{eqnarray}

(\romannumeral8) By
$$(H_{\alpha,i},H_{\beta,j},L_{\gamma,k})=(H_{\beta,j},H_{\alpha,i},L_{\gamma,k}),$$
we have
$$(H_{\alpha,i}H_{\beta,j})L_{\gamma,k}-H_{\alpha,i}(H_{\beta,j}L_{\gamma,k})=(H_{\beta,j}H_{\alpha,i})L_{\gamma,k}-H_{\beta,j}(H_{\alpha,i}L_{\gamma,k}).$$
Hence
\begin{eqnarray}\nonumber
&&d(\alpha,i,\beta,j)L_{\alpha+\beta,i+j}L_{\gamma,k}+e(\alpha,i,\beta,j)H_{\alpha+\beta,i+j}L_{\gamma,k}-c(\beta,j,\gamma,k)H_{\alpha,i}H_{\beta+\gamma,j+k}\\ \nonumber
&=&d(\beta,j,\alpha,i)L_{\alpha+\beta,i+j}L_{\gamma,k}+e(\beta,j,\alpha,i)H_{\alpha+\beta,i+j}L_{\gamma,k}-c(\alpha,i,\gamma,k)H_{\beta,j}H_{\alpha+\gamma,i+k}.
\end{eqnarray}
Then we have
\begin{eqnarray}\nonumber
&&a(\alpha+\beta,i+j,\gamma,k)(d(\alpha,i,\beta,j)-d(\beta,j,\alpha,i))\\ \nonumber
&=&c(\beta,j,\gamma,k)d(\alpha,i,\beta+\gamma,j+k)-c(\alpha,i,\gamma,k)d(\beta,j,\alpha+\gamma,j+k),\\ \nonumber
&&c(\alpha+\beta,i+j,\gamma,k)(e(\alpha,i,\beta,j)-e(\beta,j,\alpha,i))\\ \nonumber
&=&c(\beta,j,\gamma,k)e(\alpha,i,\beta+\gamma,j+k)-c(\alpha,i,\gamma,k)e(\beta,j,\alpha+\gamma,j+k).
\end{eqnarray}
By (\romannumeral3), we can conclude that
\begin{eqnarray}\nonumber
&&c(\beta,j,\gamma,k)d(\alpha,i,\beta+\gamma,j+k)-c(\alpha,i,\gamma,k)d(\beta,j,\alpha+\gamma,i+k)=0,\\ \nonumber
&&c(\beta,j,\gamma,k)e(\alpha,i,\beta+\gamma,j+k)-c(\alpha,i,\gamma,k)e(\beta,j,\alpha+\gamma,i+k)=0.
\end{eqnarray}

(\romannumeral9) By
$$(H_{\alpha,i},H_{\beta,j},H_{\gamma,k})=(H_{\beta,j},H_{\alpha,i},H_{\gamma,k}),$$
we have
$$(H_{\alpha,i}H_{\beta,j})H_{\gamma,k}-H_{\alpha,i}(H_{\beta,j}H_{\gamma,k})=(H_{\beta,j}H_{\alpha,i})H_{\gamma,k}-H_{\beta,j}(H_{\alpha,i}H _{\gamma,k}).$$
Hence
\begin{eqnarray}\nonumber
&&d(\alpha,i,\beta,j)L_{\alpha+\beta,i+j}H_{\gamma,k}+e(\alpha,i,\beta,j)H_{\alpha+\beta,i+j}H_{\gamma,k}\\ \nonumber
&&-d(\beta,j,\gamma,k)H_{\alpha,i}L_{\beta+\gamma,j+k}-e(\beta,j,\gamma,k)H_{\alpha,i}H_{\beta+\gamma,j+k}\\\nonumber
&=&d(\beta,j,\alpha,i)L_{\alpha+\beta,i+j}H_{\gamma,k}+e(\beta,j,\alpha,i)H_{\alpha+\beta,i+j}H_{\gamma,k}\\ \nonumber
&&-d(\alpha,i,\gamma,k)H_{\beta,j}L_{\alpha+\gamma,i+k}-e(\alpha,i,\gamma,k)H_{\beta,j}H_{\alpha+\gamma,i+k}.
\end{eqnarray}
Then we have
\begin{eqnarray}\nonumber
&&b(\alpha+\beta,i+j,\gamma,k)(d(\alpha,i,\beta,j)-d(\beta,j,\alpha,i))\\ \nonumber
&&+e(\alpha+\beta,i+j,\gamma,k)(e(\alpha,i,\beta,j)-e(\beta,j,\alpha,i))\\ \nonumber
&=&e(\beta,j,\gamma,k)e(\alpha,i,\beta+\gamma,j+k)+d(\alpha,i,\gamma,k)c(\beta,j,\alpha+\gamma,j+k)\\ \nonumber
&&-d(\beta,j,\gamma,k)c(\alpha,i,\beta+\gamma,j+k)-e(\alpha,i,\gamma,k)e(\beta,j,\alpha+\gamma,j+k)\\ \nonumber
&&d(\alpha+\beta,i+j,\gamma,k)(e(\alpha,i,\beta,j)-e(\beta,j,\alpha,i))\\ \nonumber
&=&e(\beta,j,\gamma,k)d(\alpha,i,\beta+\gamma,j+k)-e(\alpha,i,\gamma,k)d(\beta,j,\alpha+\gamma,j+k).
\end{eqnarray}
By (\romannumeral3), we can conclude that
\begin{eqnarray}\nonumber
&&e(\beta,j,\gamma,k)e(\alpha,i,\beta+\gamma,j+k)+d(\alpha,i,\gamma,k)c(\beta,j,\alpha+\gamma,i+k)\\ \nonumber
&=&d(\beta,j,\gamma,k)c(\alpha,i,\beta+\gamma,j+k)+e(\alpha,i,\gamma,k)e(\beta,j,\alpha+\gamma,i+k),\\ \nonumber
&&e(\beta,j,\gamma,k)d(\alpha,i,\beta+\gamma,j+k)-e(\alpha,i,\gamma,k)d(\beta,j,\alpha+\gamma,i+k)\\ \nonumber
&=&0.
\end{eqnarray}
\end{proof}
\end{lemma}

By Lemma 12, in order to obtain a compatible left-symmetric algebra structure on $\mathcal{L}(\Gamma)$, we  suppose that
\begin{eqnarray}
a(\alpha,i,\beta,j)=\frac{-\beta(1+\epsilon\beta)}{1+\epsilon(\alpha+\beta)}.
\end{eqnarray}

\begin{lemma}
Let $\mathcal{L}(\Gamma)$ be the generalized loop Heisenberg-Virasoro algebra. For $a(\alpha,i,\beta,j)$ defined by (3.4), only the following solutions satisfy all formulas in (3.3) simultaneously:
\begin{eqnarray}\nonumber
&&b(\alpha,i,\beta,j)=-\beta(1+(1-\epsilon\beta)m\delta_{\alpha+\beta,0}),\\ \nonumber
&&c(\alpha,i,\beta,j)=-\beta(1+\epsilon\beta)m\delta_{\alpha+\beta,0},\\ \nonumber
&&d(\alpha,i,\beta,j)=0,\\ \nonumber
&&e(\alpha,i,\beta,j)=0,
\end{eqnarray}
where $\alpha,\beta\in\Gamma,\,i,j\in\mathbb{Z}$ and $m$ are the fixed constant in $\mathbb{C}$.

\begin{proof}
For any $\alpha,\beta\in\Gamma,\,i,j\in\mathbb{Z}$, we  suppose that
\begin{eqnarray}
&&B(\alpha,i,\beta,j)=\frac{1+\epsilon(\alpha+\beta)}{1+\epsilon\beta}b(\alpha,i,\beta,j),\\ 
&&C(\beta,j,\alpha,i)=\frac{1+\epsilon(\alpha+\beta)}{1+\epsilon\alpha}c(\beta,j,\alpha,i).
\end{eqnarray}
Then by (3.3), we have:
\begin{eqnarray}
&&(1+\epsilon\beta)B(\alpha,i,\beta,j)-(1+\epsilon\alpha)C(\beta,j,\alpha,i)=-\beta(1+\epsilon(\alpha+\beta)),\\ \nonumber
&&B(\beta,j,\gamma,k)B(\alpha,i,\beta+\gamma,j+k)-B(\alpha,i,\gamma,k)B(\beta,j,\alpha+\gamma,i+k)\\
&&=(\alpha-\beta)B(\alpha+\beta,i+j,\gamma,k),\\ \nonumber
&&C(\beta,j,\gamma,k)B(\alpha,i,\beta+\gamma,j+k)+\gamma C(\beta,j,\alpha+\gamma,i+k)\\
&&=-\beta C(\alpha+\beta,i+j,\gamma,k).
\end{eqnarray}
Taking $\alpha=i=\beta=j=\gamma=k=0$ in (3.7) and (3.9), then 
$$B(0,0,0,0)=C(0,0,0,0),\,C(0,0,0,0)B(0,0,0,0)=0.$$
Hence $$B(0,0,0,0)=C(0,0,0,0)=0.$$
Taking $\alpha=i=0$ in (3.7) and (3.9), then
\begin{eqnarray}
&&(1+\epsilon\beta)B(0,0,\beta,j)-C(\beta,j,0,0)=-\beta(1+\epsilon\beta),\\ 
&&C(\beta,j,\gamma,k)B(0,0,\beta+\gamma,j+k)+\gamma C(\beta,j,\gamma,k)=-\beta C(\beta,j,\gamma,k).
\end{eqnarray}
Taking $\beta=\beta+\gamma,\,j=j+k$ in (3.10), then
$$B(0,0,\beta+\gamma,j+k)=\frac{C(\beta+\gamma,j+k,0,0)}{1+\epsilon(\beta+\gamma)}-(\beta+\gamma).$$
Then it is plugged into the (3.11), we obtain
$$\frac{C(\beta,j,\gamma,k)C(\beta+\gamma,j+k,0,0)}{1+\epsilon(\beta+\gamma)}-(\beta+\gamma)C(\beta,j,\gamma,k)+\gamma C(\beta,j,\gamma,k)=-\beta C(\beta,j,\gamma,k).$$
Therefore $C(\beta,j,\gamma,k)C(\beta+\gamma,j+k,0,0)=0$. 
Taking $\gamma=k=0$, then $C^2(\beta,j,0,0)=0$, thus $C(\beta,j,0,0)=0,\,B(0,0,\beta,j)=-\beta$. 
Taking $\beta=j=0,\,\alpha=-\gamma,i=-k$ in (3.9), then
$$C(0,0,\gamma,k)B(-\gamma,-k,\gamma,k)+\gamma C(0,0,0,0)=0,$$
Thus $C(0,0,\gamma,k)B(-\gamma,-k,\gamma,k)=0.$

(\romannumeral1) $B(-\alpha,-i,\alpha,i)=0,\,\mbox{for all}\,\alpha\in\Gamma,\,i\in\mathbb{Z}.$

Taking $\alpha+\beta+\gamma=0,\,i+j+k=0$ in (3.9), we have
$$\gamma C(\beta,j,-\beta,-j)=-\beta C(-\gamma,-k,\gamma,k),$$
taking the formula in (3.7), we have
$$\frac{\beta\gamma}{1-\epsilon\beta}=\frac{\beta\gamma}{1+\epsilon\gamma}.$$
This conclusion contradicts $\alpha+\beta+\gamma=0$.

(\romannumeral2) $C(0,0,\alpha,i)=0,\,\alpha\in\Gamma,\,i\in\mathbb{Z}.$

Taking $\alpha=-\gamma,\,\beta=0,\,i=-k,\,j=0$ in (3.7), we have
$$B(-\gamma,-k,0,0)=0.$$
Taking $\alpha=\beta=-\gamma,\,i=j=-k$ in (3.9), we have
$$C(-\gamma,-k,\gamma,k)B(-\gamma,-k,0,0)+\gamma C(-\gamma,-k,0,0)=\gamma C(-2\gamma,-2k,\gamma,k).$$
Since $B(-\gamma,-k,0,0)=0,\,C(\alpha,i,0,0)=0$, then
$$C(-2\gamma,-2k,\gamma,k)=0.$$
Taking $2\alpha=\beta=-2\gamma,\,2i=j=-2k$ in (3.9), we have
$$C(-2\gamma,-2k,\gamma,k)B(-\gamma,-k,-\gamma,-k)+\gamma C(-2\gamma,-2k,0,0)=2\gamma C(-3\gamma,-3k,\gamma,k).$$
Since $C(-2\gamma,-2k,\gamma,k)=0,\,C(\alpha,i,0,0)=0$, then
$$C(-3\gamma,-3k,\gamma,k)=0.$$
\ldots\ldots\\
Taking $a\alpha=\beta=-a\gamma,\,ai=j=-ak,\,a\in\mathbb{C}$ in (3.9), we have
\begin{eqnarray}\nonumber
&&C(-a\gamma,-ak,\gamma,k)B(-\gamma,-k,(-a+1)\gamma,(-a+1)k)+\gamma C(-a\gamma,-ak,0,0)\\ \nonumber
&=&a\gamma C(-(-a+1)\gamma,-(-a+1)k,\gamma,k).
\end{eqnarray}
Then $$C(-(-a+1)\gamma,-(-a+1)k,\gamma,k)=0.$$
We can conclude that $$C(-a\alpha,-ai,\alpha,i)=0,\,a\in\mathbb{C}\setminus\{1\}$$
By (3.7), we have
$$B(\alpha,i,-a\alpha,-ai)=a\alpha\frac{1+\epsilon(\alpha-a\alpha)}{1-\epsilon a\alpha},\,a\in\mathbb{C},\,a\geq2.$$
Taking $\alpha=0,\,\beta=-\gamma$ in (3.8), we have 
$$B(-\gamma,-k,\gamma,k)=0.$$
Similarly, taking $-\alpha=\beta=\gamma,\,-i=j=k$ in (3.9), we have
$$C(\gamma,k,\gamma,k)B(-\gamma,-k,2\gamma,3k)+\gamma C(\gamma,k,0,0)=-\gamma C(\gamma,k,\gamma,k).$$
Since $C(0,0,\gamma,k)=0,\,C(\gamma,k,0,0)=0,\,B(-\gamma,-k,2\gamma,2k)\neq0$, then
$$C(\gamma,k,\gamma,k)=0.$$
Taking $-2\alpha=\beta=2\gamma,\,-2i=j=2k$ in (3.9), we have
$$C(2\gamma,2k,\gamma,k)B(-\gamma,-k,3\gamma,3k)+\gamma C(2\gamma,2k,0,0)=-2\gamma C(\gamma,k,\gamma,k).$$
Since $C(\gamma,k,\gamma,k)=0=0,\,C(\gamma,k,0,0)=0,\,B(-\gamma,-k,3\gamma,3k)\neq0$, then
$$C(2\gamma,k,2\gamma,k)=0.$$
\ldots\ldots\\
Taking $-a\alpha=\beta=a\gamma,\,-ai=j=ak,\,a\in\mathbb{C}$ in (3.9), we have
$$C(a\gamma,ak,\gamma,k)B(-\gamma,-k,(a+1)\gamma,(a+1)k)+\gamma C(a\gamma,ak,0,0)=-a\gamma C((a-1)\gamma,(a-1)k,\gamma,k).$$
Then $$C(a\gamma,ak,\gamma,k)=0.$$
We can conclude that $$C(a\alpha,ai,\alpha,i)=0,\,a\in\mathbb{C}\setminus\{-1\}.$$
Thus
$$C(\alpha,i,\beta,j)=0,\,\alpha,\beta\in\Gamma,\,i,j\in\mathbb{Z},\,\alpha+\beta\neq0.$$
Taking $\alpha+\beta+\gamma=0,\,i+j+k=0$ in (3.9), since $B(-\gamma,-k,\gamma,k)=0$, then we have
$$\gamma C(\beta,j,-\beta,-j)=-\beta C(-\gamma,-k,\gamma,k).$$
Taking $\beta=j=1$, we have $C(-\gamma,-k,\gamma,k)=-\gamma C(1,1,-1,-1).$
Let $m=C(1,1,-1,-1)$. Then  
$$C(-\gamma,-k,\gamma,k)=-\gamma m.$$
When $\alpha+\beta=0$, then $C(\alpha,i,\beta,j)=-\beta m$, taking it in (3.6), we have 
$$c(\alpha,i,\beta,j)=-\beta m(1+\epsilon\beta),\,\alpha+\beta=0.$$
Taking in (3.3), we have
$$b(\alpha,i,\beta,j)=-\beta(1+(1-\epsilon\beta)m),\,\alpha+\beta=0.$$
When $\alpha+\beta\neq0$, then $C(\alpha,i,\beta,j)=0$, taking it in (3.6), we have
$$c(\alpha,i,\beta,j)=0,\,\alpha+\beta\neq0.$$
Taking in (3.3), we have 
$$b(\alpha,i,\beta,j)=-\beta,\,\alpha+\beta\neq0.$$
Thus
$$b(\alpha,i,\beta,j)=-\beta(1+(1-\epsilon\beta)m\delta_{\alpha+\beta,0}),\,c(\alpha,i,\beta,j)=-\beta(1+\epsilon\beta)m\delta_{\alpha+\beta,0}.$$
Then, we calculate $d(\alpha,i,\beta,j),\,e(\alpha,i,\beta,j)$. Taking $\alpha=i=0$ in seventh formula of (3.3), we have
$$-(\beta+\gamma)d(\beta,j,\gamma,k)=-\beta d(\beta,j,\gamma,k).$$
We have $$d(\beta,j,\gamma,k)=0.$$
Similarly, taking $\alpha=i=0$ in eighth formula of (3.3), we have 
$$-(\beta+\gamma)e(\beta,j,\gamma,k)=-\beta e(\beta,j,\gamma,k).$$
We have $$e(\beta,j,\gamma,k)=0.$$
Conversely, it is obvious.
\end{proof} 
\end{lemma}

\end{document}